\documentclass[12pt]{article}
\usepackage{amssymb,amsmath}
\usepackage{hyperref}
\usepackage{graphicx}
\usepackage{color}
\usepackage{chngcntr}
\usepackage{upquote}

\begin{document}

\title{\LARGE\bf A sampling-based approximation \\
of the complex error function and its implementation without poles}

\author{
\normalsize\bf S. M. Abrarov\footnote{\scriptsize{Dept. Earth and Space Science and Engineering, York University, Toronto, Canada, M3J 1P3.}}\,, B. M. Quine$^{*}$\footnote{\scriptsize{Dept. Physics and Astronomy, York University, Toronto, Canada, M3J 1P3.}}\, and R. K. Jagpal$^{\dagger}$}

\date{April 3, 2018}
\maketitle

\begin{abstract}
Recently we developed a new sampling methodology based on incomplete cosine expansion of the sinc function and applied it in numerical integration in order to obtain a rational approximation for the complex error function $w\left(z \right) = e^{- {z^2}}\left(1 + \frac{2i}{\sqrt \pi}\int_0^z e^{t^2}dt\right),$ where $z = x + iy$. As a further development, in this work we show how this sampling-based rational approximation can be transformed into alternative form for efficient computation of the complex error function $w\left(z \right)$ at smaller values of the imaginary argument $y=\operatorname{Im}\left[z \right]$. Such an approach enables us to avoid poles in implementation and to cover the entire complex plain with high accuracy in a rapid algorithm. An optimized Matlab code utilizing only three rapid approximations is presented.
\vspace{0.25cm}
\vspace{0.25cm}
\end{abstract}

\section{Introduction}

The complex error function, also known as the Faddeeva function, can be defined as \cite{Faddeyeva1961, Abramowitz1972}
\begin{equation}\label{eq_1}
w\left( z \right) = {e^{ - {z^2}}}\left( {1 + \frac{{2i}}{{\sqrt \pi  }}\int\limits_0^z {{e^{{t^2}}}dt} } \right),
\end{equation}
where $z = x + iy$. Despite simple representation, the integral in equation \eqref{eq_1} cannot be taken analytically. Therefore, the integral equation \eqref{eq_1} for the complex error function must be computed numerically.

This function $w\left( z \right)$ is the most important in family of the Faddeeva functions. In particular, comparing this equation with the error function of complex argument \cite{Abramowitz1972}
$$
{\text{erf}}\left( z \right) = \frac{2}{{\sqrt \pi  }}\int\limits_0^z {{e^{{-t^2}}}dt}
$$
one can see that
$$
{\text{erf}}\left( z \right) = 1 - {e^{ - {z^2}}}w\left( {iz} \right).
$$
Consequently, we can conclude that the error function ${\text{erf}}\left( z \right)$ is just a reformulation of the complex error function $w\left( z \right)$. Some other functions of practical importance that can be reformulated in terms of the complex error function are the Dawson's integral \cite{Cody1970, McCabe1974, Rybicki1989, Nijimbere2017, Abrarov2018}

$$
{\text{daw}}\left( z \right) = {e^{ - {z^2}}}\int\limits_0^z {{e^{{t^2}}}dt}  = \frac{{\sqrt \pi  }}{{2i}}\left( {w\left( z \right) - {e^{ - {z^2}}}} \right),
$$
the Fresnel integral \cite{Abramowitz1972}
$$
F\left( z \right) = \int\limits_0^z {{e^{i\left( {\pi /2} \right){t^2}}}dt = \left( {1 + i} \right)\left[ {1 - {e^{i\left( {\pi /2} \right){z^2}}}w\left( {\sqrt \pi  \left( {1 + i} \right)z/2} \right)} \right]} /2
$$
and the plasma dispersion function \cite{Fried1961}
$$
\mathcal{Z}\left( z \right) = PV\frac{1}{{\sqrt \pi  }}\int\limits_{ - \infty }^\infty  {\frac{{{e^{ - {t^2}}}}}{{t - z}}dt}  = i\sqrt \pi  w\left( z \right),
$$
where the principal value signifies that it has no discontinuity at $y = \operatorname{Im} \left[ z \right] = 0$.

The equation \eqref{eq_1} is absolutely identical to the complex probability function on the upper half of the complex plain \cite{Armstrong1972} (see also \cite{Abramowitz1972})
$$
W\left( z \right) = PV\frac{i}{\pi }\int\limits_{ - \infty }^\infty  {\frac{{{e^{ - {t^2}}}}}{{z - t}}dt}  = w\left( z \right), \qquad \operatorname{Im} \left[ z \right] \geqslant 0,
$$
where the principal value implies again that this function has no discontinuity at $y = \operatorname{Im} \left[ z \right] = 0$. The real part of the complex probability function is known as the Voigt function \cite{Armstrong1967, Schreier1992, Letchworth2007, Pagnini2010, Abrarov2015a}
$$
K\left( {x,y} \right) = \operatorname{Re} \left[ {W\left( z \right)} \right] = PV\frac{y}{\pi }\int\limits_{ - \infty }^\infty  {\frac{{{e^{ - {t^2}}}}}{{{y^2} + {{\left( {x - t} \right)}^2}}}dt},
$$
widely used in atmospheric science to describe the spectral line broadening effects in photon absorption and emission by various gas molecular species in a planetary atmosphere.

There is a remarkable property of the complex error function \cite{Abramowitz1972, Gautschi1970, McKenna1984}
\begin{equation}\label{eq_2}
w\left( z \right) = 2{e^{ - {z^2}}} - w\left( { - z} \right).
\end{equation}
From equation \eqref{eq_2} it immediately follows that
$$
w\left( {x, - \left| y \right|} \right) = 2{e^{ - {{\left( { - x + i\left| y \right|} \right)}^2}}} - w\left( { - x,\left| y \right|} \right)
$$
This signifies that it is sufficient to consider only I$^\text{st}$ and II$^\text{nd}$ quadrants in order to cover the entire complex plane. Thus, in order to simplify algorithmic implementation we will imply further that $y \geqslant 0$.

The identity \eqref{eq_2} can be used not only to avoid a direct computation of the complex error function at negative $y < 0$. In this work we show how this identity can also be generalized to derive a sampling-based approximation of the complex error function that excludes all its poles in algorithmic implementation. This approach sustains high accuracy in computation at smaller values of the parameter $y$ that is commonly considered difficult for computation of the Voigt/complex error function \cite{Armstrong1967, Amamou2013, Abrarov2015b}.

\section{Methodology and derivation}

In our earlier publication \cite{Abrarov2015c} we have developed a new methodology of sampling based on incomplete cosine expansion of the sinc function. This technique of sampling is especially efficient in numerical integration. As an example, we have shown that applying the incomplete cosine expansion of the sinc function to equation \eqref{eq_1}, the following sampling-based rational approximation of the complex error function
\begin{equation}\label{eq_3}
w\left( z \right) \approx \sum\limits_{m = 1}^{{2^{K - 1}}} {\frac{{{A_m} + {B_m}\left( {z + i\varsigma /2} \right)}}{{C_m^2 - {{\left( {z + i\varsigma /2} \right)}^2}}}},
\end{equation}
where the expansion coefficients are
$$
{A_m} = \frac{{\sqrt \pi  \left( {2m - 1} \right)}}{{{2^{2K}}h}}\sum\limits_{n =  - N}^N {{e^{{\varsigma ^2}/4 - {n^2}{h^2}}}\sin \left( {\frac{{\pi \left( {2m - 1} \right)\left( {nh + \varsigma /2} \right)}}{{{2^K}h}}} \right)},
$$
$$
{B_m} =  - \frac{i}{{{2^{K - 1}}\sqrt \pi  }}\sum\limits_{n =  - N}^N {{e^{{\varsigma ^2}/4 - {n^2}{h^2}}}\cos \left( {\frac{{\pi \left( {2m - 1} \right)\left( {nh + \varsigma /2} \right)}}{{{2^K}h}}} \right)}
$$
and
$$
{C_m} = \frac{{\pi \left( {2m - 1} \right)}}{{{2^{K + 1}}h}},
$$
can be obtained. Specifically, by taking $h = 0.25$, $\varsigma  = 2.75$, $K = 5$ and $N = 23$ the approximation \eqref{eq_3} alone covers with average accuracy $\sim 10^{-14}$ \cite{Abrarov2015c} the entire domain $0 \leqslant x \leqslant 40,000$ and ${10^{ - 4}} \leqslant y \leqslant {10^2}$ of the Voigt line-shapes that arise from the collection of HITRAN molecular transitions \cite{Rothman2013} at low terrestrial altitudes and for spectral displacements less than $25 \,\, cm^{-1}$.

Later we have shown that the truncation integer in approximation \eqref{eq_3} may not be necessarily equal to ${2^{K - 1}}$. In particular, this restriction can be avoided by replacing ${2^{K - 1}}$ with an arbitrary integer $M$ (see \cite{Abrarov2015d} for more details). Consequently, the approximation \eqref{eq_3} can be rewritten in form
\begin{equation}\label{eq_4}
w\left( z \right) \approx \sum\limits_{m = 1}^M {\frac{{{a_m} + {b_m}\left( {z + i\varsigma /2} \right)}}{{c_m^2 - {{\left( {z + i\varsigma /2} \right)}^2}}}},
\end{equation}
where the expansion coefficient are modified correspondingly as
$$
{a_m} = \frac{{\sqrt \pi  \left( {m - 1/2} \right)}}{{2{M^2}h}}\sum\limits_{n =  - N}^N {{e^{{\varsigma ^2}/4 - {n^2}{h^2}}}\sin \left( {\frac{{\pi \left( {m - 1/2} \right)\left( {nh + \varsigma /2} \right)}}{{Mh}}} \right)},
$$
$$
{b_m} =  - \frac{i}{{M\sqrt \pi  }}\sum\limits_{n =  - N}^N {{e^{{\varsigma ^2}/4 - {n^2}{h^2}}}\cos \left( {\frac{{\pi \left( {m - 1/2} \right)\left( {nh + \varsigma /2} \right)}}{{Mh}}} \right)} 
$$
and
$$
{c_m} = \frac{{\pi \left( {m - 1/2} \right)}}{{2Mh}}.
$$

Although the approximation \eqref{eq_4} can cover the entire HITRAN domain, its accuracy deteriorates as the parameter $y$ decreases. In order to resolve this problem we rewrite equation \eqref{eq_2} as follows
\begin{equation}\label{eq_5}
w\left( z \right) = {e^{ - {z^2}}} + \frac{{w\left( z \right) - w\left( { - z} \right)}}{2}.
\end{equation}
Substituting approximation \eqref{eq_4} into right side of the equation \eqref{eq_5}, after some trivial rearrangements we can transform it into alternative form as follows
\begin{equation}\label{eq_6}
w\left( z \right) \approx e^{-z^2} + z\sum\limits_{m = 1}^M {\frac{{{\alpha _m} - {\beta _m}{z^2}}}{{{\gamma _m} - {\theta _m}{z^2} + {z^4}}}},
\end{equation}
where
$$
{\alpha _m} = {b_m}\left[c_m^2 - \left(\frac{\varsigma ^2}{2}\right)^2 \right] + i{a_m}\varsigma = {b_m}\left[ {{{\left( {\frac{{\pi \left( {m - 1/2} \right)}}{{2Mh}}} \right)}^2} - {{\left( {\frac{\varsigma }{2}} \right)}^2}} \right] + i{a_m}\varsigma,
$$
$$
{\beta _m} = {b_m},
$$
$$
{\gamma _m} = c_m^4 + \frac{c_m^2{\varsigma ^2}}{2} + \frac{\varsigma ^4}{16} = {\left[ {{{\left( {\frac{{\pi \left( {m - 1/2} \right)}}{{2Mh}}} \right)}^2} + {{\left( {\frac{\varsigma }{2}} \right)}^2}} \right]^2}
$$
and
$$
{\theta _m} = 2c_m^2 - \frac{\varsigma ^2}{2} = 2{\left( {\frac{{\pi \left( {m - 1/2} \right)}}{{2Mh}}} \right)^2} - {\frac{\varsigma }{2}^2}.
$$
Since the new equation \eqref{eq_6} is derived by transformation from the sampling-based rational approximation \eqref{eq_4}, it also represents a sampling-based approximation.

It should be noted that application of the identity \eqref{eq_5} in derivation of approximation of the complex error function in alternative form has been proposed already in our recent publication \cite{Abrarov2016} (see also the corresponding Matlab code \cite{Matlab47801}). However, in this work we show its generalization leading to approximation \eqref{eq_6} that can be used without poles in a rapid algorithm. Furthermore, we also suggest that equation \eqref{eq_5} may be applied to other approximations in order to compute more accurately the Voigt/complex error function at smaller values of the parameter $y$.

\section{Implementation}

\subsection{Approximations and boundaries}

Similar to our previous work \cite{Abrarov2018} we applied only three approximations bounded inside domains as shown in Fig. 1; in fact, due to symmetric properties of the complex error function
$$
\operatorname{Re} \left[ {w\left( {x, y} \right)} \right] = \operatorname{Re} \left[ {w\left( { - x, y} \right)} \right]
$$
and
$$
\operatorname{Im} \left[ {w\left( {x, y} \right)} \right] =  - \operatorname{Im} \left[ {w\left( { - x, y} \right)} \right],
$$
it is sufficient to consider only the I$^\text{st}$ quadrant of the complex plain. The complex plain is divided into external and internal domains. Taking as an objective a worst relative accuracy of ${10^{ - 13}}$, we separated these domains by boundaries accordingly for the best optimization.

Internal domain is situated inside a circle $\left| {x + iy} \right| = 8$ and consists of two subdomains, the primary subdomain and secondary subdomain. The secondary subdomain is bounded by a straight line $y = 0.05\left| x \right|$.

It is very convenient for algorithmic implementation to rewrite the approximation \eqref{eq_4} as follows
\begin{equation}\label{eq_7}
\begin{aligned}
 w\left( z \right) \approx & \,\Omega \left( {z + i\varsigma /2} \right) \\ 
 \Rightarrow & \,\Omega \left( z \right) \triangleq \sum\limits_{m = 1}^M {\frac{{{a_m} + {b_m}z}}{{c_m^2 - {z^2}}}}. 
\end{aligned}
\end{equation}
At $M = 23$ this approximation meets the requirement for accuracy exceeding ${10^{ - 13}}$ within primary subdomain.

For secondary subdomain we may apply the approximation \eqref{eq_6} without any modification. However, at $M = 23$ its accuracy becomes $\sim{10^{ - 12}}$ in the area near the origin. In order to resolve this problem we should increase the number of summation terms by two as given by \footnote{These two additional terms may be optional if the requirement for accuracy $\sim 10^{-12}$ is sufficient for users.}
\begin{equation}\label{eq_8}
w\left( z \right) \approx e^{-z^2} + z\sum\limits_{m = 1}^{M + 2} {\frac{{{\alpha _m} - {\beta _m}{z^2}}}{{{\gamma _m} - {\theta _m}{z^2} + {z^4}}}}.
\end{equation}
Inclusion of these two terms almost does not decelerate the computation and sustains high accuracy exceeding ${10^{ - 13}}$ everywhere within the secondary domain.

External domain utilizes the following Laplace continued fraction \cite{Abramowitz1972, Gautschi1970, Poppe1990a} given by
\begin{equation}\label{eq_9}
w\left( z \right) \approx \frac{{\left( {i/\sqrt \pi  } \right)}}{{z - \frac{{1/2}}{{z - \frac{1}{{z - \frac{{3/2}}{{z - \frac{2}{{z - \frac{{5/2}}{{z - \frac{3}{{z - \frac{{7/2}}{{z - \frac{4}{{z - \frac{{9/2}}{{z - \frac{5}{{z - \frac{{11/2}}{z}}}}}}}}}}}}}}}}}}}}}}}.
\end{equation}

An optimized Matlab code, implemented according to this scheme, is shown in Appendix A. \\

\begin{figure}[ht]
\begin{center}
\includegraphics[width=20pc]{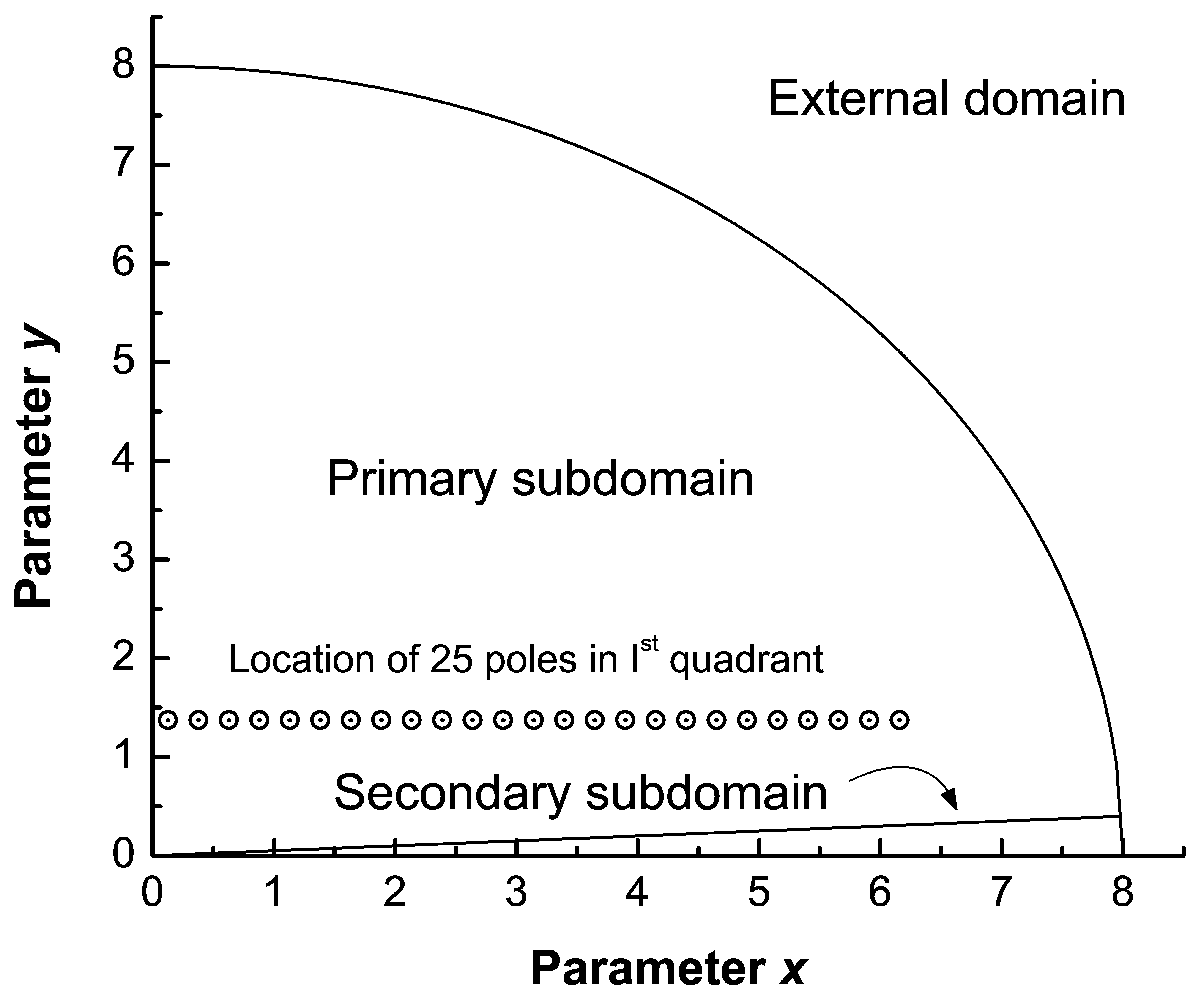}\hspace{1pc}%
\begin{minipage}[b]{28pc}
\vspace{0.3cm}
{\sffamily {\bf{Fig. 1.}} Boundaries and location of $25$ poles in ${\text{I}}^{\text{st}}$ quadrant of the complex plain.}
\end{minipage}
\end{center}
\end{figure}

\vspace{-0.75cm}

\subsection{Location of poles}

Approximation \eqref{eq_7} contains poles that can be readily found by solving the following quadratic equation
$$
c_m^2 - {\left( {z + i\varsigma /2} \right)^2} = 0.
$$
It is not difficult to see that the solution for the equation above results to two poles at each index $m$. Particularly, we can find that
$$
{z_{1,2}} =  \pm c_m - i\varsigma/2.
$$
Fortunately, all these poles are located in the ${\text{II}}{{\text{I}}^{{\text{rd}}}}$ and ${\text{I}}{{\text{V}}^{{\text{th}}}}$ quadrants only since $\operatorname{Im} \left[ {{z_{1,2}}} \right] =  - \varsigma /2 =  - 1.375$. Therefore, these poles do not affect the computation.

The new sampling-based approximation \eqref{eq_8} also contains poles. In particular, solving the following quartic equation
\begin{equation}
\label{eq_10}
{\gamma _m} - {\theta _m}{z^2} + {z^4} = 0
\end{equation}
one can find four poles associated with each index $m$. These poles are
$$
{z_{1,2}} =  \pm \frac{{\sqrt {{\theta _m} + \sqrt {\theta _m^2 - 4{\gamma _m}} } }}{{\sqrt 2 }}
$$
and
$$
{z_{3,4}} =  \pm \frac{{\sqrt {{\theta _m} - \sqrt {\theta _m^2 - 4{\gamma _m}} } }}{{\sqrt 2 }}.
$$
All these poles are located along two horizontal lines, since (see Appendix B)
$$
\operatorname{Im} \left[ {{z_{1,4}}} \right] = \varsigma /2 = 1.375
$$
and
$$
\operatorname{Im} \left[ {{z_{2,3}}} \right] = - \varsigma /2 = - 1.375.
$$

The location of $M+2=25$ poles on the ${{\text{I}}^{{\text{st}}}}$ quadrant is shown in Fig. 1 by open circles with dots inside. Since these poles are situated far beyond the secondary subdomain, they also do not affect the computation.

\subsection{Error analysis}

The error analysis is performed by using relative errors defined as
$$
{\Delta _{\operatorname{Re} }} = \left| {\frac{{\operatorname{Re} \left[ {{w_{ref}}\left( {x,y} \right)} \right] - \operatorname{Re} \left[ {w\left( {x,y} \right)} \right]}}{{\operatorname{Re} \left[ {{w_{ref}}\left( {x,y} \right)} \right]}}} \right|
$$
and
$$
{\Delta _{\operatorname{Im} }} = \left| {\frac{{\operatorname{Im} \left[ {{w_{ref}}\left( {x,y} \right)} \right] - \operatorname{Im} \left[ {w\left( {x,y} \right)} \right]}}{{\operatorname{Im} \left[ {{w_{ref}}\left( {x,y} \right)} \right]}}} \right|,
$$
where ${w_{ref}}\left( {x,y} \right)$ is the reference, for the real and imaginary parts, respectively.

\begin{figure}[ht]
\begin{center}
\includegraphics[width=32pc]{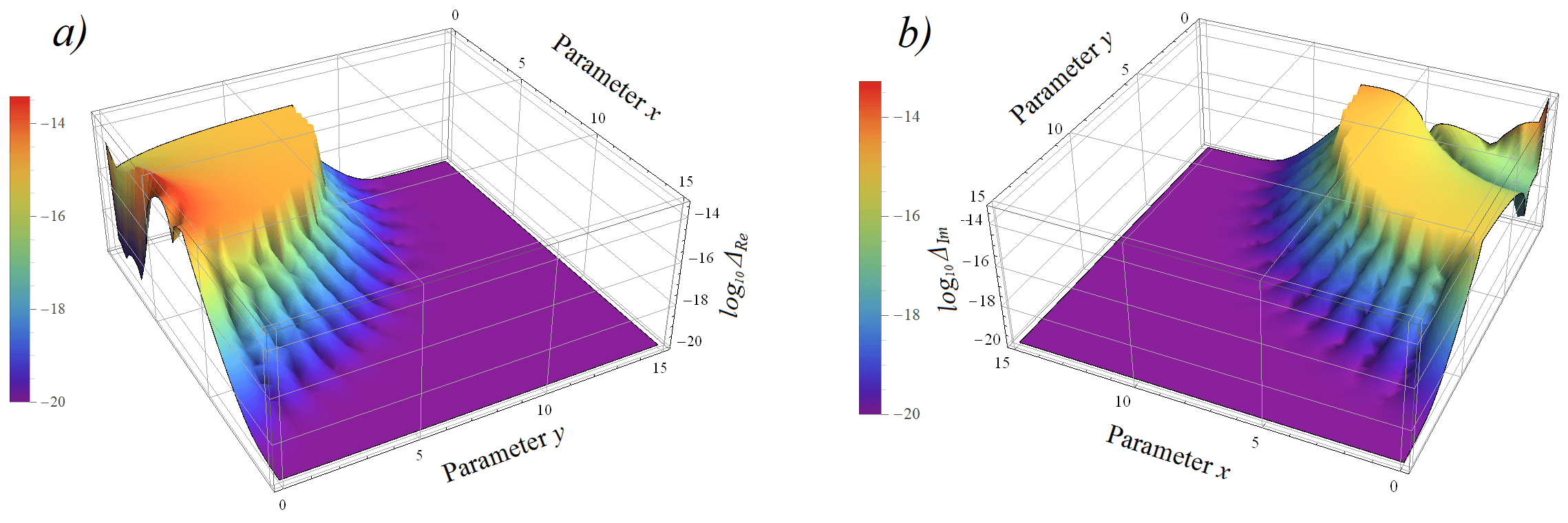}\hspace{1pc}%
\begin{minipage}[b]{28pc}
\vspace{0.3cm}
{\sffamily {\bf{Fig. 2.}} Logarithm of the relative errors for (a) real and (b) imaginary parts over the domain $0\leqslant x \leqslant 15$ and $0 \leqslant y \leqslant 15$.}
\end{minipage}
\end{center}
\end{figure}

\bigskip
\begin{figure}[ht]
\begin{center}
\includegraphics[width=32pc]{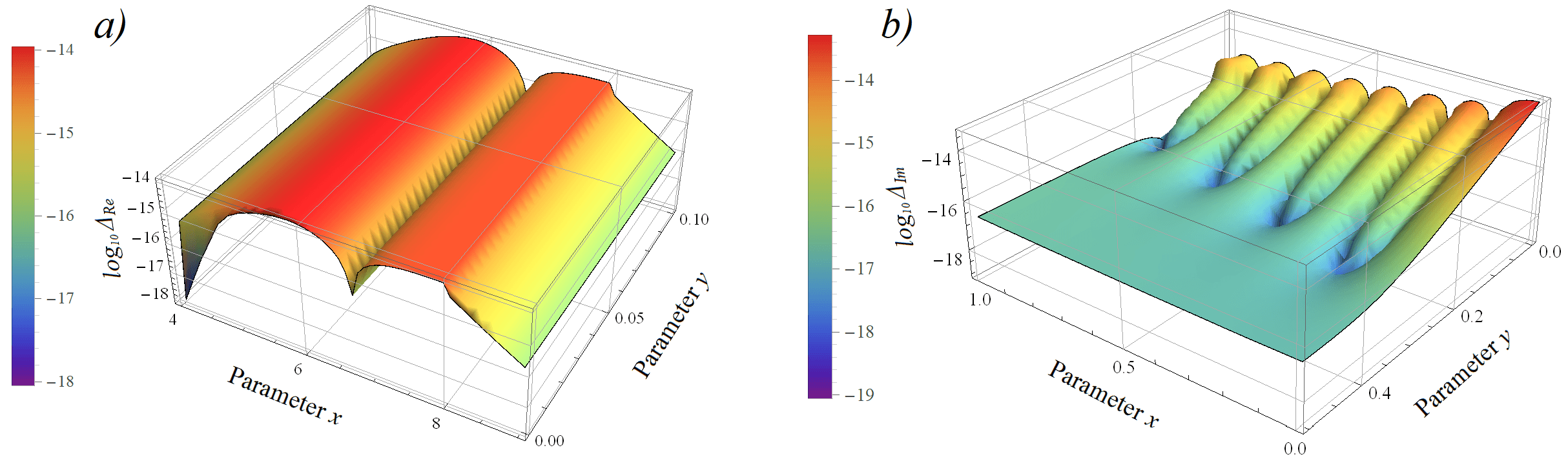}\hspace{1pc}%
\begin{minipage}[b]{28pc}
\vspace{0.3cm}
{\sffamily {\bf{Fig. 3.}} Logarithm of the relative errors over the domains (a) $4\leqslant x \leqslant 10$ and $0\leqslant y \leqslant 0.1$, (b) $0\leqslant x \leqslant 1$ and $0\leqslant y \leqslant 0.5$ for the real and imaginary parts, respectively.}
\end{minipage}
\end{center}
\end{figure}

Figures 2a and 2b show the logarithm of relative errors for the real and imaginary parts, respectively, over the area $0 \leqslant x \leqslant 15$ and $0 \leqslant y \leqslant 15$. Figures 3a and 3b depict the logarithm of relative errors with worst accuracies. As we can see, over the real part the worst accuracy is about $2 \times {10^{ - 14}}$ while over the imaginary part the worst accuracy is about $8 \times {10^{ - 14}}$.

Thus, the provided accuracy at double precision computation of the complex error function $w \left( z \right)$ is absolutely consistent with {\it{CERNLIB}}, {\it{libcerf}} and {\it{RooFit}} packages (see \cite{Karbach2014} for detailed information regarding accuracy of these packages).

\section{Run-time test}

The Matlab is an array programing language. Consequently, the number of applied equations should be minimized in order to reduce computational flow consisting of large size arrays \cite{Abrarov2018}. The algorithm we develop employs only three rapid approximations \eqref{eq_7}, \eqref{eq_8} and \eqref{eq_9}. Although equation \eqref{eq_8} involves quartic polynomial in its denominator, this practically does not decelerate the computation since the arrays ${z^2}$, ${z^4}$ can be predefined prior to the nested sum (see the function {\ttfamily{SD2 = subdom2(z)}} inside body of the Matlab code in Appendix A). Furthermore, the array of the exponential function $e^{-z^2}$ is computed just once outside the nested loop. Consequently, the approximation \eqref{eq_8} is almost as rapid as the rational approximation \eqref{eq_7}.

The run-time test has been performed by comparing our Matlab code, shown in Appendix A, with Ab-Initio group of MIT implementation written in C/C++ programing language by Steven Johnson \cite{Johnson, Matlab38787}. The corresponding algorithm in this C/C++ implementation is based on some modification of the Algorithm $680$ \cite{Poppe1990a, Poppe1990b} with additional inclusion of the Salzer's approximation for smaller values of the imaginary part $y = \operatorname{Im} \left[ z \right]$ (see our recent paper \cite{Abrarov2018} for detailed information, see also \cite{Salzer1951, Zaker1969}).

Although Ab-Initio group of MIT implementation utilizes only earliest equations published for the complex error function, it can provide, nevertheless, a rapid computation with relative errors smaller than ${10^{ - 13}}$ \cite{Johnson}.

Despite that Matlab programs are usually slower than their C/C++ analogs, the run-time test shows that with $10$ million random numbers, generated within the domain $0 < x < 6$ and $0 < y < 0.1$, our Matlab code is faster by a factor $1.68$ ($11.6$ and $6.9$ seconds, respectively). This can be explained from the fact that the Salzer's equation covering this domain in C/C++ implementation is not simple and, therefore, relatively slow due to requirement to compute multiple times the hyperbolic sine and hyperbolic cosine functions in a nested recurrence procedure (see our recent work \cite{Abrarov2018} for C/C++ and Matlab command lines and technical aspects describing why this domain is especially important for consideration).

Computational test also reveals that with $10$ million random numbers, generated within the domain $\left| {x + iy} \right| < 15$, the Matlab code is faster than C/C++ code by a factor $1.51$ ($5.9$ and $3.9$ seconds, respectively). Within the domain $\left| {x + iy} \right| < 10,000$ the C/C++ implementation is faster than the Matlab code by a factor $1.33$ ($2.1$ and $2.8$ seconds, respectively, for $10$ million random numbers). However, since a high spectral resolution in computational spectroscopy is required only for the sharp lines and absolutely unnecessary for the flat curve regions, this advantage of the C/C++ implementation disappears when we perform a computation with non-equidistantly distributed grid-points. For example, if we take, say, $9$ million grid points over the domain $\left| {x + iy} \right| < 15$ and the remaining $1$ million grid points beyond this domain, then the Matlab code remains faster than the C/C++ implementation \cite{Johnson, Matlab38787}. It should be noted that application of the non-equidistantly distributed grid-points is very common in radiative transfer applications since this approach significantly accelerates the computation \cite{Smith1978, Quine2002}. In atmospheric radiative transfer applications the Loretzian distribution can be used to select non-equidistantly spaced grid-points along $x$-axis at some fixed value $y$.

The present Matlab code is as fast as that of reported in our recent publication \cite{Abrarov2018} with some very minor differences over specific domains. To the best of our knowledge this and our recent \cite{Abrarov2018} Matlab codes are the most rapid in computation of the complex error function as compared to any other Matlab codes with comparable accuracy ever reported in scientific literature or elsewhere online (see for example this \cite{Matlab12091} and this \cite{Matlab47801} Matlab codes). Therefore, the proposed approximations \eqref{eq_7} and \eqref{eq_8} may be also be useful in radiative transfer models like MODRTAN \cite{Berk2017} and bytran \cite{Pliutau2017}, where our previously reported approximations \cite{Abrarov2011, Abrarov2015d} are currently used for rapid and accurate computation of the Voigt/complex error function.

The run-time test has been performed on a typical desktop computer (Intel(R), CPU at $2.6$ GHz, RAM $8$ GB, Windows $10$) by using the Intel compiler and Matlab $2009$b distributive.

\section{Conclusion}

In this work we show how the sampling-based approximation can be obtained and applied for efficient computation of the complex error function $w\left( z \right)$ at smaller values of the imaginary argument $y = \operatorname{Im} \left[ z \right]$. This approach results in coverage of the entire complex plain. An optimized Matlab code providing rapid computation with high accuracy is presented.

\section*{Acknowledgments}

This work is supported by National Research Council Canada, Thoth Technology Inc. and York University.

\vspace{-0.25cm}
\section*{Appendix A}

function FF = fadsamp(z)

\footnotesize
\begin{verbatim}
%     This function file computes the complex error function (also known as
% the Faddeeva function) by using a new method of sampling based on
% incomplete cosine expansion of the sinc function [1, 2]. External domain 
% is computed by the Laplace continued fraction [3]. The description of the
% algorithm is presented in the work [4].
%
% REFERENCES
% [1] S. M. Abrarov and B. M. Quine, Sampling by incomplete cosine
%     expansion of the sinc function: application to the Voigt/complex
%     error function, Appl. Math. Comput., 258 (2015) 425-435.
%     https://doi.org/10.1016/j.amc.2015.01.072
%
% [2] S. M. Abrarov and B. M. Quine, A rational approximation for efficient
%     computation of the Voigt function in quantitative spectroscopy, J.
%     Math. Research, 7 (2) (2015) 163-174.
%     https://doi.org/10.5539/jmr.v7n2p163
%
% [3] W. Gautschi, Efficient computation of the complex error function. SIAM
%     J. Numer. Anal., 7 (1) (1970) 187-198.
%     https://doi.org/10.1137/0707012
%
% [4] S. M. Abrarov, B. M. Quine and R. K. Jagpal, A sampling-based
%     approximation of the complex error function and its implementation
%     without poles, Appl. Numer. Math., 129 (2018) 181-191.
%     https://doi.org/10.1016/j.apnum.2018.03.009
%
%     The code is written by Sanjar M. Abrarov, Brendan M. Quine and
% Rajinder K. Jagpal, York University, Canada, February 2018.
%
% *************************************************************************
% All variables in this section are global within the function file.
% *************************************************************************

h = 0.25; % step
stigma = 2.75; % shift constant
m_max = 23; % truncating integer for index m
n_max = 23; % truncating integer for index n

n = -n_max:n_max; % array for index n

tab = ones(6,m_max); % initiate the table
m = 1; % counter
while m <= m_max + 2 % the expansion coefficients

    tab(1,m) = (sqrt(pi)*(m - 1/2))/(2*m_max^2*h)*sum(exp(stigma^2/4 - ...
        n.^2*h^2).*sin((pi*(m - 1/2).*(n*h + stigma/2))/(m_max*h)));

    tab(2,m) = -1i/(m_max*sqrt(pi))*sum(exp(stigma^2/4 - ...
        n.^2*h^2).*cos((pi*(m - 1/2).*(n*h + stigma/2))/(m_max*h)));

    tab(3,m) = pi*(m - 1/2)/(2*m_max*h);

    tab(4,m) = tab(2,m)*(((pi*(m - 1/2))/(2*m_max*h))^2 - ...
        (stigma/2)^2) + 1i*tab(1,m)*stigma;

    tab(5,m) = (((pi*(m - 1/2))/(2*m_max*h))^2 + (stigma/2)^2)^2;

    tab(6,m) = 2*((pi*(m - 1/2))/(2*m_max*h))^2 - stigma^2/2;

    m = m + 1; % increment the counter

    % The expansion coefficients are:
    % a = tab(1,:);
    % b = tab(2,:);
    % c = tab(3,:);
    % alpha = tab(4,:);
    % beta  = tab(2,:);
    % gamma = tab(5,:);
    % theta = tab(6,:);
end
% *************************************************************************
% End of section.
% *************************************************************************

ind_neg = imag(z)<0; % if some imag(z) values are negative, then ...
z(ind_neg) = conj(z(ind_neg)); % ... bring them to the upper-half plane

FF = zeros(size(z));

ind_ext  = abs(z)>8; % external indices

FF(~ind_ext) = intf(z(~ind_ext));
FF(ind_ext) = contfr(z(ind_ext)); % continued fraction (external region)

function CF = contfr(z) % the Laplace continued fraction approximation
    
    bN = 11; % initial integer
    bN = 1:bN;
    bN = bN/2;

    CF = bN(end)./z; % start computing from the last bN
    for k = 1:length(bN) - 1
        CF = bN(end-k)./(z - CF);
    end
    CF = 1i/sqrt(pi)./(z - CF);
end

function IF = intf(z) % internal function

IF = zeros(size(z));
    
ind_pr = imag(z)>0.05*abs(real(z));

IF(ind_pr) = Omega(z(ind_pr) + 1i*stigma/2);
IF(~ind_pr) = subdom2(z(~ind_pr)); % secondary subdomain
    
    function OF = Omega(z) % Omega function for primary subdomain
        
        zP2 = z.^2; % define repeating array

        OF = 0;
        for k = 1:n_max
            OF = OF + (tab(1,k) + tab(2,k)*z)./(tab(3,k).^2 - zP2);
        end
    end

    function SD2 = subdom2(z) % secondary subdomain

        zP2 = z.^2; % first repeating array
        zP4 = zP2.^2; % second repeating array

        SD2 = 0;
        for k = 1:n_max + 2 % increased by 2 terms!
            SD2 = SD2 + (tab(4,k) - tab(2,k)*zP2)./(tab(5,k) - tab(6,k)* ...
                zP2 + zP4);
        end
        SD2 = exp(-z.^2) + z.*SD2;
    end
end

% Convert for negative imag(z) values
FF(ind_neg) = conj(2*exp(-z(ind_neg).^2) - FF(ind_neg));
end
\end{verbatim}
\normalsize

\section*{Appendix B}

There are four solutions of the quartic equation \eqref{eq_10}. In order to find them it is convenient to represent equation \eqref{eq_10} in a biquadratic form as follows
\[
\label{eq_B.1}
\tag{B.1}
{\gamma _m} - {\theta _m}Z + {Z^2} = 0,
\]
where
\[
\label{eq_B.2}
\tag{B.2}
Z = {z^2}.
\]

Two solutions of the biquadratic equation \eqref{eq_B.1} are given by
\[
\label{eq_B.3}
\tag{B.3}
{Z_{1,2}} = \frac{{{\theta _m} \pm \sqrt {\theta _m^2 - 4{\gamma _m}} }}{2}.
\]

Since
$$
{\theta _m} = 2c_m^2 - \frac{{{\varsigma ^2}}}{2}
$$
and
$$
{\gamma _m} = c_m^4 + \frac{{c_m^2{\varsigma ^2}}}{2} + \frac{{{\varsigma ^4}}}{{16}}
$$
it follows that
$$
\theta _m^2 - 4{\gamma _m} =  - 4{c^2}{\varsigma ^2}.
$$
Consequently, from equation \eqref{eq_B.3} we have
\small
\[
\label{eq_B.4}
\tag{B.4}
{Z_{1,2}} = \frac{{\left( {2c_m^2 - \varsigma^2 /2} \right) \pm \sqrt { - 4c_m^2{\varsigma ^2}} }}{2} = \frac{{\left( {2c_m^2 - \varsigma^2 /2} \right) \pm 2i{c_m}\varsigma }}{2} = \frac{1}{4}{\left( {2{c_m} \pm i\varsigma } \right)^2}.
\]
\normalsize
Lastly, taking into account the relation \eqref{eq_B.2} from equation \eqref{eq_B.4} we obtain four solutions for the equation \eqref{eq_10}
$$
{z_1} = \frac{1}{2}\left( {2{c_m} + i\varsigma } \right) = c_m + i\frac{\varsigma}{2},
$$
$$
{z_2} = - \frac{1}{2}\left( {2{c_m} + i\varsigma } \right) = - c_m - i\frac{\varsigma}{2},
$$
$$
{z_3} = \frac{1}{2}\left( {2{c_m} - i\varsigma } \right) = c_m - i\frac{\varsigma}{2},
$$
and
$$
{z_4} = - \frac{1}{2}\left( {2{c_m} - i\varsigma } \right) = - c_m + i\frac{\varsigma}{2}.
$$

\bigskip


\begin{thebibliography}{9}

\bibitem{Faddeyeva1961}
V.N. Faddeyeva, and N.M. Terent\text{'}ev, Tables of the probability integral $w\left( z \right) = {e^{ - {z^2}}}\left( {1 + \frac{{2i}}{{\sqrt \pi  }}\int_0^z {{e^{{t^2}}}dt} } \right)$ for complex argument. Pergamon Press, Oxford, 1961.

\bibitem{Abramowitz1972}
M. Abramowitz and I.A. Stegun. Error function and Fresnel integrals. Handbook of mathematical functions with formulas, graphs, and mathematical tables. $9^{\text{th}}$ ed. New York 1972, 297-309.

\bibitem{Cody1970}
W.J. Cody, K.A. Paciorek and H.C. Thacher, Chebyshev approximations for Dawson\text{'}s integral. Math. Comp. 24 (1970) 171-178. \\
\url{https://doi.org/10.1090/S0025-5718-1970-0258236-8}

\bibitem{McCabe1974}
J.H. McCabe, A continued fraction expansion with a truncation error estimate for Dawson\text{'}s integral, Math. Comp. 28 (1974) 811-816. \\ 
\url{https://doi.org/10.1090/S0025-5718-1974-0371020-3}

\bibitem{Rybicki1989}
G.B. Rybicki, Dawson\text{'}s integral and the sampling theorem, Comp. Phys., 3 (1989) 85-87. \\
\url{https://doi.org/10.1063/1.4822832}

\bibitem{Nijimbere2017}
V. Nijimbere, Analytical evaluation and asymptotic evaluation of Dawson's integral and related functions in mathematical
physics, \href{https://arxiv.org/abs/1703.06757}{arXiv:1703.06757} (2017).

\bibitem{Abrarov2018}
S.M. Abrarov and B.M. Quine, A rational approximation of the Dawson's integral for efficient computation of the complex error function
Applied Math. Comput., 321 (15) (2018) 526-543. \\
\url{https://doi.org/10.1016/j.amc.2017.10.032}

\bibitem{Fried1961}
B.D. Fried and S.D. Conte. The plasma dispersion function. New York, Academic Press, 1961.

\bibitem{Armstrong1972}
B.H. Armstrong and B.W. Nicholls, Emission, absorption and transfer of radiation in heated atmospheres. Pergamon Press, New York, 1972.

\bibitem{Armstrong1967}
B.H. Armstrong, Spectrum line profiles: The Voigt function, J. Quantit. Spectrosc. Radiat. Transfer, 7 (1) (1967) 61-88. \\
\url{https://doi.org/10.1016/0022-4073(67)90057-X}

\bibitem{Schreier1992}
F. Schreier, The Voigt and complex error function: A comparison of computational methods. J. Quant. Spectrosc. Radiat. Transfer, 48 (5-6) (1992) 743-762. \\
\url{https://doi.org/10.1016/0022-4073(92)90139-U}

\bibitem{Letchworth2007}
K.L. Letchworth and D.C. Benner, Rapid and accurate calculation of the Voigt function, J. Quant. Spectrosc. Radiat. Transfer, 107 (1) (2007) 173-192. \\
\url{https://doi.org/10.1016/j.jqsrt.2007.01.052}

\bibitem{Pagnini2010}
G. Pagnini and F. Mainardi, Evolution equations for the probabilistic generalization of the Voigt profile function, J. Comput. Appl. Math., 233 (6) (2010) 1590-1595. \\
\url{https://doi.org/10.1016/j.cam.2008.04.040}

\bibitem{Abrarov2015a}
S.M. Abrarov and B.M. Quine, A rational approximation for efficient computation of the Voigt function in quantitative spectroscopy, J. Math. Research, 7 (2) (2015) 163-174. \\
\url{https://doi.org/10.5539/jmr.v7n2p163}

\bibitem{Gautschi1970}
W. Gautschi, Efficient computation of the complex error function. SIAM J. Numer. Anal., 7 (1) (1970) 187-198. \\
\url{https://doi.org/10.1137/0707012}

\bibitem{McKenna1984}
S.J. McKenna, A method of computing the complex probability function and other related functions over the whole complex plane. Astrophys. Space Sci., 107 (1) (1984) 71-83. \\
\url{https://doi.org/10.1007/BF00649615}

\bibitem{Amamou2013}
H. Amamou, B. Ferhat and A. Bois, Calculation of the Voigt function in the region of very small values of the parameter $a$ where the calculation is notoriously difficult, Amer. J. Anal. Chem., 4 (2013) 725-731. \\
\url{https://doi.org/10.4236/ajac.2013.412087}

\bibitem{Abrarov2015b}
S.M. Abrarov and B.M. Quine, Accurate approximations for the complex error function with small imaginary argument, J. Math. Research 7 (1) (2015) 44-53. \\
\url{https://doi.org/10.5539/jmr.v1n1p44}

\bibitem{Abrarov2015c}
S.M. Abrarov and B.M. Quine, Sampling by incomplete cosine expansion of the sinc function: Application to the Voigt/complex error function, Appl. Math. Comput., 258 (2015) 425-435. \\
\url{https://doi.org/10.1016/j.amc.2015.01.072}

\bibitem{Rothman2013}
L.S. Rothman, I.E. Gordon, Y. Babikov, A. Barbe, D.C. Benner, P.F. Bernath, M. Birk, L. Bizzocchi, V. Boudon, L.R. Brown, A. Campargue, K. Chance, E.A. Cohen, L.H. Coudert, V.M. Devi, B.J. Drouin, A. Fayt, J.-M. Flaud, R.R. Gamache, J.J. Harrison, J.-M. Hartmann, C. Hill, J.T. Hodges, D. Jacquemart, A. Jolly, J. Lamouroux, R.J. Le Roy, G. Li, D.A. Long, O.M. Lyulin, C.J. Mackie, S.T. Massie, S. Mikhailenko, H.S.P. M\"{u}ler, O.V. Naumenko, A.V. Nikitin, J. Orphal, V. Perevalov, A. Perrin, E.R. Polovtseva and C. Richard, The HITRAN2012 molecular spectroscopic database, J. Quant. Spectrosc. Radiat. Transfer, 130 (2013) 4-50. \\
\url{https://doi.org/10.1016/j.jqsrt.2013.07.002}

\bibitem{Abrarov2015d}
S.M. Abrarov and B.M. Quine, A rational approximation for efficient computation of the Voigt function in quantitative spectroscopy,
J. Math. Research, 7 (2) (2015) 163-174. \\
\url{https://doi.org/10.5539/jmr.v7n2p163}

\bibitem{Abrarov2016}
S.M. Abrarov and B.M. Quine, The Fourier expansion approximation for high-accuracy computation of the Voigt/complex error function at small imaginary argument, \href{https://arxiv.org/abs/1606.07871}{arXiv:1606.07871} (2016).

\bibitem{Matlab47801}
\href{https://www.mathworks.com/matlabcentral/fileexchange/47801-the-voigt-complex-error-function--second-version-}{Matlab Central, file ID \#: 47801} (2016).

\bibitem{Poppe1990a}
G.P.M. Poppe and C.M.J. Wijers, More efficient computation of the complex error function. ACM Transact. Math. Software, 16 (1990) 38-46. \\
 \url{https://doi.org/10.1145/77626.77629}

\bibitem{Karbach2014}
T.M. Karbach, G. Raven and M. Schiller, Decay time integrals in neutral meson mixing and their efficient evaluation, \href{https://arxiv.org/abs/1407.0748}{arXiv:1407.0748} (2014).

\bibitem{Johnson}
S. G. Johnson, Faddeeva package (2017). \\
\url{http://ab-initio.mit.edu/wiki/index.php/Faddeeva_Package}

\bibitem{Matlab38787}
\href{https://www.mathworks.com/matlabcentral/fileexchange/38787-faddeeva-package--complex-error-functions}{Matlab Central, file ID \#: 38787} (2012).

\bibitem{Poppe1990b}
G.P.M. Poppe and C.M.J. Wijers, Algorithm 680: evaluation of the complex error function. ACM Transact. Math. Software, 16 (1990) 47. \\
\url{https://doi.org/10.1145/77626.77630}
 
\bibitem{Salzer1951}
H.E. Salzer, Formulas for calculating the error function of a complex variable, Math. Tables Aids Comput. 5 (1951) 61-70. \\
\url{https://doi.org/10.2307/2002163}

\bibitem{Zaker1969}
T.A. Zaker, Calculation of the complementary error function of complex argument, J. Comput. Phys., 4 (3) (1969) 427-430. \\
\url{https://doi.org/10.1016/0021-9991(69)90011-4}

\bibitem{Smith1978}
H.J.P. Smith, D.J. Dube, M.E. Gardner, S.A. Clough and F.X. Kneizys, FASCODE - Fast Atmospheric Signature Code (Spectral Transmittance and Radiance), \href{http://www.dtic.mil/dtic/tr/fulltext/u2/a057506.pdf}{AFGL-TR-78-0081, AD A057506} (1978).

\bibitem{Quine2002}
B.M. Quine and J.R. Drummond, GENSPECT: a line-by-line code with selectable interpolation error tolerance J. Quant. Spectrosc. Radiat. Transfer 74 (2002) 147-165. \\
\url{https://doi.org/10.1016/S0022-4073(01)00193-5}

\bibitem{Matlab12091}
\href{https://www.mathworks.com/matlabcentral/fileexchange/12091-complex-scaled-complementary-error-function}{Matlab Central, file ID \#: 12091} (2007).

\bibitem{Berk2017}
A. Berk and F. Hawes, Validation of MODTRAN{\textregistered}6 and its line-by-line algorithm, J. Quantit. Spectrosc. Radiat. Transfer, 203 (2017) 542-556.\\
\url{https://doi.org/10.1016/j.jqsrt.2017.03.004}

\bibitem{Pliutau2017}
D. Pliutau and K. Roslyakov, Bytran $-{\mid}-$ spectral calculations for portable devices using the HITRAN database, Earth Sci. Informatics, 2017 (10) (3) 395-404. \\
\url{https://doi.org/10.1007/s12145-017-0288-4}

\bibitem{Abrarov2011}
S.M. Abrarov and B.M. Quine, Efficient algorithmic implementation of the Voigt/complex error function based on exponential series approximation, Appl. Math. Comput., 218 (5) (2011) 1894-1902. \\
\url{https://doi.org/10.1016/j.amc.2011.06.072}

\end{thebibliography}
\end{document}